\documentclass[11pt]{ip-journal}

\usepackage{amsmath, amsthm, amsfonts, verbatim, amssymb}
\usepackage{tikz}
\usetikzlibrary{matrix,arrows}
\usepackage{url}
\usepackage[all,cmtip]{xy}
\usepackage[hidelinks]{hyperref}
\hypersetup{
     colorlinks   = true,
     citecolor    = blue
}
\hypersetup{linkcolor=blue}
\usepackage{xcolor}
\usepackage{graphicx}
\usepackage{pinlabel}

\newcommand\bN{{\mathbb N}}

\def\cM{{\mathcal M}}
\def\ocM{\overline{\mathcal M}}

\def\cO{{\mathcal O}}
\def\cU{\mathcal U}

\def\cL{{\mathcal L}}

\def\cX{\mathcal X}
\def\bR{{\mathbb R}}
\def\bZ{{\mathbb Z}}

\def\bC{{\mathbb C}}

\def\cM{{\mathcal M}}

\def\cC{{\mathcal C}}

\newtheorem {theo}{Theorem}
\newtheorem {coro}{Corollary}
\newtheorem {lemm}{Lemma}

\newtheorem {defi}{Definition}

\newtheorem {prop}{Proposition}

\def\tC{\widetilde C}

\def\bs{\bigskip}
\def\ms{\medskip}
\def\ni{\noindent}

\def\pd{\partial}

\def\Pic{\mbox{Pic}}

\def\oz1{d{\overline z}^1}
\def\oz2{d{\overline z}^2}
\def\oz3{d{\overline z}^3}

\def\oI{\overline I}

\def\oz{\overline z}

\def\oIq1{\oI_1\cdots\oI_{q-1}}
\def\oIq2{\oI_1\cdots\oI_{q-2}}

\def\dim{{\mbox{dim}}}

\def\hp{\hat{p}}

\def\hX{\widehat{X}}

\def\span{\mbox{span}}
\def\Av{\mbox{Av}}
\def\Res{\mbox{Res}}
\def\ord{\mbox{ord}}
\def\cV{\mathcal V}

\begin{document}

\ni
\title[Limit of Weierstrass measures]
{Limit of Weierstrass measure on stable curves}

\author[Ngai-Fung Ng and Sai-Kee Yeung]
{Ngai-Fung Ng and Sai-Kee Yeung}

\begin{abstract}  
{\it The goal of the paper is to study the limiting behavior of the Weierstrass measures on a smooth curve of genus $g\geqslant 2$ as the curve approaches a certain nodal stable curve represented by a point in the Deligne-Mumford compactification $\ocM_g$ of the moduli $\cM_g$, including irreducible ones
or those of compact type. As a consequence, the Weierstrass measures on a stable rational 
curve at the boundary of $\cM_g$ are completely determined. In the process, the asymptotic behavior of the Bergman measure is also studied.}
\end{abstract}

\address[]{N.-F. Ng, Department of Mathematics, Purdue University, West Lafayette, IN 47907 USA}
\address[]{S.-K. Yeung, Department of Mathematics, Purdue University, West Lafayette, IN 47907 USA}
\email{ngn@alumni.purdue.edu}
\email{yeung@math.purdue.edu}

\thanks{
\noindent{
The second author was partially supported by a grant from the National Science Foundation.
}}

\ni{\it}
\maketitle

\begin{center}
{\bf 1. Introduction} 
\end{center} 

\ms
\ni{\bf 1.1}  On a compact Riemann surface, an interesting geometric object to study is
the distribution of Weierstrass points associated to
the tensor powers of an ample line bundle.  It is observed by Olsen \cite{olsen} that the asymptotic distribution of such Weierstrass points is dense with respect to the analytic topology.  The situation is clarified by the 
beautiful result of Mumford \cite{mumford1} 
and Neeman \cite{neeman} that the asymptotic distribution
is uniformly distributed with respect to the Bergman kernel of the curve.  In other words, the higher Weierstrass points as defined are weakly equidistributed with respect to the Arakelov measure.  The phenomenon is interesting both from  a geometric and an arithmetic point of view, such as results
explained in \cite{burnol}, \cite{delcentina}, \cite{mumford1} and \cite{rudnick}.

A natural problem is what happens for the 
corresponding distribution on a singular algebraic curve,
in particular for stable curves at the boundary of a moduli space in its Deligne-Mumford compactification.
It has been observed by \cite{ballico-gatto}, \cite{garcia-lax}, \cite{lax}, \cite{little-furio} that
the asymptotic distribution of the Weierstrass points
associated to the power of an ample line bundle is
no longer dense with respect to the complex topology on 
a rational nodal curve.
The main goal of this paper is to clarify the situation and
give a precise statement about the distribution of the 
Weierstrass points on a nodal curve sitting at the boundary of a moduli space.

Let $\cM_g$ be the moduli space of compact Riemann surfaces of genus $g\geqslant 2$.  Let $\ocM_g$ be the
Deligne-Mumford compactification of $\cM_g$.  A point $t\in \cM_g$ represents a Riemann surface $X_t$ of genus $g$.
Consider now a one-parameter family of Riemann surfaces $\pi:\cX\rightarrow S$ on a curve
$S\subset \ocM$ so that $o\in \ocM-\cM$ and a neighborhood 
$U$ of $o$ satisfies $U-\{o\}\subset \cM$.  Our goal is to study the behavior of the limit of the Weierstrass measure on $X_o$.
Our approach is to
study the relation between the Bergman kernel and the Weierstrass measure for stable curves arising from degeneration of a family of
smooth curve.   For this purpose, we have to study 
the limiting behavior of the geometry of the period mapping and extract from it the geometric information needed.

\ms
\ni{\bf 1.2}   We refer the reader to Section 2 for various terminology
used in the introduction.  For our statement here,
the Bergman measure $\mu^B_X$ on a compact Riemann surface $X$ is defined by
$\mu^B_X=\sqrt{-1}\sum_{i=1}^g\omega_i\wedge\overline{\omega_i}$ over an orthonormal basis $\{\omega_1,\dots,\omega_g\}$ of $\Gamma(X,K_X)$.
Our first result is the estimates on the asymptotic behavior of
the Bergman measure.

\begin{theo} \label{bergmanlimit}
(a) Let $\pi:\cX\rightarrow S$ be a local family of stable curves in the sense of Deligne-Mumford so that $X_t=\pi^{-1}(t)$ is smooth for $t\in S-\{o\}$ and $X_o$ has a single
node at $p\in X_o$.
Let $\tau:\hX_o\rightarrow X_o$ be the normalization of $X_o$.  Then\\
(i) if the node $p$ on $X_o$ is separating,  the Bergman measure $\mu^B_{X_t} \rightarrow \tau_*\mu^B_{\hX_o}$ as $t\rightarrow o$.\\
(ii) if the node $p$ on $X_o$ is non-separating,  the Bergman measure $\mu^B_{X_t} \rightarrow \tau_*\mu^B_{\hX_o}+\delta_p$ as $t\rightarrow o$,
where $\delta_p$ is the Dirac Delta at the node $p\in X_o$.\\
(b) Let $X_o$ be a stable curve for which $p_1,\cdots, p_k$ are non-separating nodes and $p_{k+1},\cdots,p_l$ are
separating nodal points. Assume that $X_o$ has $l-k+1$ irreducible components. Let $\hX_o$ be the normalization of $X_o$.  Then in terms of the notations above,
$$\mu^B_{X_t}\rightarrow \mu^B_{X}=\tau_*\mu^B_{\hX}+\sum_{i=1}^k\delta_{p_i}$$
as $t\rightarrow o$
\end{theo}

In the above, we denote by $\tau_*\mu^B_{\hX_o}$ the measure 
$(\tau^{-1})^*|_{\tau^{-1}(X_o-\{p\})} \mu^B_{\hX_o}$, using the fact that
$\tau^{-1}$ is a biholomorphism on
$X_o-\{p\}$.\\

Here we remark that for a stable curve $X_o$ given by a point at the boundary of
the Deligne-Mumford compactification $D:=\pd\cM_g=\ocM_g-\cM_g$ 
of $\cM_g$ studied in this paper, the limit $\mu^B_{X_t} \rightarrow \tau_*\mu^B_{\hX_o}$ is independent of
the family of smooth curves taken.

\ms
\ni{\bf 1.3}  
 Let $L$ be an invertible sheaf of positive degree on a stable curve $X_o$ of genus $g\geqslant 2$.  The notion of Weierstrass points of powers of $L$ has been generalized from smooth curves to stable curves in the literature, cf. \cite{widland}.
Assume $X_o$ is represented by a point $o$ in the boundary of 
the Deligne-Mumford compactification of $\cM_g$.   Assume that $L$ could be extended as an invertible sheaf $L_t$ to each curve $X_t$ represented by a point $t$ in a neighborhood of $U$ of $o$ in $\ocM_g$.  It is at this juncture that
we need to assume that $X_o$ is irreducible or of compact type for a general line bundle.  In general it may be difficult to extend a line bundle $L$ consistently to a nearby fiber due to the difficulty of defining limit linear series on stable curves.  This is however possible in the case that $X_o$ is irreducible or
of compact type, cf. \cite{ak}, \cite{ce}, \cite{cp}, \cite{gz}.  A typical example is given by tensor power of the (relative) dualizing sheaf of the family.  The second author is grateful to Samuel Grushevsky for pointing out the subtlety of extension of the line bundle.

From the work of \cite{neeman} and \cite{mumford1}, if $X_t$ is a smooth curve, the discrete measure $\mu^W_{m L_{t}}$ associated to the set of Weierstrass points of $(L_{X_t})^m\to X_t$ converges to $1/g\cdot \mu^B_{X_t}$ as $m \to \infty$ where $\mu^B_{X_t}$ is the Bergman measure of $X_t$.

\begin{theo} \label{weierstrasslimit}
Let $\pi:\cX\rightarrow S$ be a local family of stable curves which are either irreducible or of compact type so that $X_t=\pi^{-1}(t)$ is smooth for $t\in S-\{o\}$ and $X_o$ has a single
node at $p\in X_o$.
Let $\hX_o$ be the normalization of $X_o$.  Then\\
(a) if the node $p$ on $X_o$ is separating, the measure $\mu^W_{mL_o}$ associated to the Weierstrass points on $X_o$ satisfies $\mu^W_{mL_o} \to 1/g \cdot \mu^B_{\hX_o} $ as $m \to \infty$.\\
(b) if the node $p$ on $X_o$ is non-separating, the measure $\mu^W_{mL_o}$ associated to the Weierstrass points on $X_o$ satisfies $\mu^W_{mL_o} \to 1/g \cdot (\mu^B_{\hX_o}+\delta_p )$ as $m \to \infty$.
\end{theo}

\ms
\ni{\bf 1.5}  The following result is a consequence of Theorems \ref{bergmanlimit}, \ref{weierstrasslimit} and induction.

\begin{theo} \label{mainthm} 
(a) Let $X$ be a stable curve which is irreducible and $p_1,\cdots, p_k$ are the non-separating nodes. Let $\hX$ be the normalization of $X$.  Then
$$\mu_{X,L} := \lim_{m\to\infty}\mu^W_{mL_X}= 1/g\cdot(\sum_{i=1}^k\delta_{p_i} + \mu^B_{\hX})$$ 
where $g$ is the genus of the curve from which $X$ is obtained by pinching corresponding cycles. \\
(b) Let $X$ be a stable curve of compact type and $p_1,\cdots, p_k$ are the separating nodes. Let $\hX$ be the normalization of $X$.  Then
$$\mu_{X,L}^W = 1/g\cdot\mu^B_{\hX}$$ 
where $g$ is the genus of the curve from which $X$ is obtained by pinching corresponding cycles.
\end{theo}

Note that the theorem shows that the measure is
independent of $L$.

\ms
\ni{\bf 1.6} As an immediate corollary, we have the following result in the case of a rational nodal curve living on the boundary of the Deligne-Mumford compactification of the moduli space of curves.

\begin{coro} \label{cor1}
Let $X$ be a stable irreducible rational nodal curve with nodes at $p_1,\cdots, p_g$.  Then $\mu_X = 1/g \cdot (\sum_{i=1}^g\delta_{p_i})$.
\end{coro}

Related to the corollary, we remark that from the earlier work of \cite{ballico-gatto}, \cite{garcia-lax}, \cite{lax}, and  \cite{little-furio}, 
it is known that $\mu_X$ vanishes on $X$ except possibly on a finite number of circles and the nodes. The corollary above shows that the measure $\mu_X$ is solely supported on the nodes.  Thereom {\ref{mainthm}} and Corollary \ref{cor1} above complete the picture on asymptotic distribution of Weierstrass points
for stable curves. 



\ms
\ni{\bf 1.7}  In the following we outline the main steps of 
proof.  Theorem \ref{bergmanlimit} follows from a careful study of the Bergman metric with respect to the degeneration at a single
node.  Theorem \ref{weierstrasslimit} is the main result.  It follows from the following three steps.  
The first is the convergence of the Weierstrass measure as one approaches
the boundary of the moduli. The second is to prove uniform
convergence of the Weierstrass measure to the Bergman measure on compacta in the complement of the nodes, which
depends on the results of Neeman \cite{neeman}.  Finally we apply Theorem \ref{bergmanlimit} to deduce that the residual measure is 
supported at a node.  Theorem \ref{mainthm} follows from Theorem \ref{weierstrasslimit} and an induction argument.

\ms
\ni{\bf 1.8} {\it Acknowledgement}  The authors are very grateful to Samuel Grushevsky for making many valuable comments and suggestions on the paper, to Valery Alexeev for explaining his work on semiabelic pairs. The authors would also like to thank the referee for helpful comments on the article.

\ms
\ni{\bf 1.9} After the paper was accepted, we were kindly informed of the earlier articles \cite{am}, \cite{dj} which are related to the study of the Bergman measure in general, see Remark 16.4 of \cite{dj}, and distribution of the Weierstrass point over tropical curves.  See also related more recent work in \cite{shivaprasad}.

\bs
\begin{center}
{\bf 2. Preliminaries} 
\end{center} 

\bs
\ni{\bf 2.1} Denote by $\cM_g$ the moduli space of Riemann surfaces of genus $g\geqslant 2$.  Let $\ocM_g$ be the Deligne-Mumford compactification of $\cM_g$.  The points on the boundary $\ocM_g-\cM_g$ represent stable curves in the sense of Deligne-Mumford.  For simplicity of notation, sometimes we just
denote $\cM_g, \ocM_g$ by $\cM, \ocM$ when there is no danger of confusion.

It is well-known that the compactifying 
divisor $D=\ocM_g-\cM_g$ has a decomposition 
$D=D_0\cup\cdots D_{[g/2]}$ into irreducible components, where $[x]$ denotes the integral part of 
$x$.  The generic point of the stratum $D_0$
represents an irreducible complex curve of genus $g-1$ with $2$ points identified.  We call such a node non-separating.  
The generic point of the stratum $D_i$ for $i>0$ represents the union of two irreducible complex curves of genus $j$ and $g-j$, each with a puncture
and the two punctures are identified, named as a separating node on the curve.

The nodes are obtained by contracting a real $1$-cycle on a smooth 
Riemann surface of genus $g$, by considering a family of curves $C_t$ of genus $g$ which is smooth for $t\neq 0$ and $C_0$ is the stable nodal curve considered.

A generic point on the intersection of two components $D_i\cap D_j$ for $i\neq j$ corresponds to a stable curve obtained by contracting two real cycles to two
different nodes.

We refer the readers to \cite{hm} for basic facts
about moduli space of curves.

\ms
\ni{\bf 2.2} Let $X$ be a compact Riemann surface of genus $g\geqslant 1$.  The space of holomorphic one forms
$\Gamma(X,K_X)$ has dimension $g$.  There is a natural $L^2$ metric on $\Gamma(X,K_X)$ defined by
$(\eta_1,\eta_2)=\sqrt{-1}\int_X\eta_1\wedge\overline{\eta_2}.$  We would denote by $\{\omega_1,\dots,\omega_g\}$ an orthonormal
basis of $\Gamma(X,K_X)$ on $X$.  

\begin{defi}
The Bergman measure $\mu^B_X$ on $X$ is defined by
$\mu^B_X=\sqrt{-1}\sum_{i=1}^g\omega_i\wedge\overline{\omega_i}$ where $\{\omega_1,\dots,\omega_g\}$ is any orthonormal basis of $\Gamma(X,K_X)$.
\end{defi}

It is a standard fact that the Bergman measure is independent of the orthonormal basis chosen. The Bergman measure $\mu^B_X$ is also given by the pull back of the flat measure on the Jacobian of the Riemann surface $X$ by the Abel-Jacobian map.

\ms
\ni{\bf 2.3}
A symplectic homology basis of $X$ is a basis $\{A_j, B_j\}_{1\leqslant j\leqslant g}$ of $H_1(X,\bZ)$ 
satisfying intersection pairings
$$A_i\cdot A_j=0,  \  B_i\cdot B_j=0 \ ,  \ \mbox{and}  \ 
A_i\cdot B_j=\delta_{ij} \ \mbox{for all $i$ and $j$}.$$

A canonically normalized basis $\{\omega'_i\}_{i=1}^g$ of  $\Gamma(X,K_X)$ with respect to a symplectic homology basis $\{A_j, B_j\}$ is a basis satisfying $\int_{A_j}\omega_k'=\delta_{jk}$ for all $j$ and $k$.

The period matrix of $X$ is the $g\times g$ matrix defined by
$\Omega_{ij}=\int_{B_i}\omega_j'$.

In the following we recall some standard results on the behavior of canonically normalized holomorphic one forms with respect to the symplectic bases of a deformation family. We refer the reader to \cite{fay}, \cite{yamada} and \cite{wentworth} for any unexplained terminology.

Consider a one-parameter family of Riemann surfaces $\pi:\cX\rightarrow S$ with
$S-\{o\}\subset \cM$ and  $o\in \ocM-\cM$.  The stable nodal curve $X_o$ is a singular curve with nodal points
as the only singularities, which are also called punctures of $X_o$.  $X_o$ can be considered as a union of finitely many
compact Riemann surfaces $\hX_o$ with some particular points identified corresponding to the nodal points 
where $\hX_o$ is the normalization of $X_o$.  The local defining equation for a neighborhood of a node can be
described by $zw=0$ in $\bC^2$.

We will assume that $X_o$ has only one node.  The cases of more than one node will follow from induction.

We recall the limit of canonically normalized holomorphic one forms as $t\rightarrow o$. There are two cases according to whether a node is separating or non-separating.  

\ms
\ni{\bf 2.4} In the case of separating node $p$, $X_t$ degenerates into $X_1\cup X_2$ as $t\rightarrow o$, where $X_1$ and $X_2$ are
Riemann surfaces of genus $g_1=g(X_1)>0$ and $g_2=g(X_2)>0$, and $p$ is represented by
$x_1\in X_1$ and $x_2\in X_2$. Here $g=g_1+g_2.$  Analytic structure of the degeneration is understood from the following model.  For $i=1,2$, let $x_i\in X_i$ representing $p$ and $U_i$ be a neighborhood of $x_i$ in $X_i$ with coordinates $z_i:U_i\rightarrow\Delta$
centered at $x_i$.  Let $S=\{(x,y,t):xy=t, x,y,t\in\Delta_1\}$.
Denote the fiber at $t\in \Delta$ by $S_t$. Here $\Delta_r$ denotes the disk of radius $r$ in $\bC$.  For $|t|<1$, 
glue together $X_1-z_1^{-1}(\Delta_{|t|})$ and $X_2-z_2^{-1}(\Delta_{|t|})$ according to the recipe
$$z_1 \mapsto (z_1,\frac t{z_1},t), \ z_2\mapsto (\frac t{z_2},z_2,t).$$  The resulting surfaces give rise
to an analytic family $\cX\rightarrow \Delta_1$ with smooth fibers $X_t$ for $t\neq 0$ centered around $X_0=X_o.$
For $z\in X_i-\{p\}$ and $|t|$ sufficiently small, there is
a natural section $z(t)$ of $\cX\rightarrow \Delta_1$ with $z(0)=z$.  In the notation of \cite{wentworth}, we say that
$z\in X_i\cap X_t$ if $z(t)\in X_i-z_i^{-1}(\Delta_{|t|^{1/2}})$ for all small $t$.

Let $\{\omega^{(1)\prime}_i\}$, $\{\omega^{(2)\prime}_j\}$ be normalized bases with respect to some symplectic
homology bases on $X_1$ and $X_2$ respectively.

\begin{prop}\label{sep} (\cite{wentworth} page 433, \cite{fay} page 38, \cite{yamada} page 129)
We can find a normalized basis of $\Gamma(X_t,K_{X_t})$ for $t$ sufficiently close to $o$ such that
$$
\omega'_i(x,t)=\left\{\begin{array}{lcc}
\omega^{(1)\prime}_i(x)+\cO(t^2),&\mbox{for}& x\in X_1-U_1,\\
-t\omega^{(1)\prime}_i(x)\omega^{(2)\prime}(x,p)+\cO(t^2),&\mbox{for}& x\in X_2-U_2,
\end{array}\right.
$$
$$
\omega'_j(x,t)=\left\{\begin{array}{lcc}
\omega^{(2)\prime}_j(x)+\cO(t^2),&\mbox{for}& x\in X_2-U_2,\\
-t\omega^{(2)\prime}_j(x)\omega^{(1)\prime}(x,p)+\cO(t^2),&\mbox{for}& x\in X_1-U_1,\end{array}\right.
$$
where $1\leqslant i\leqslant g_1$, $g_1+1\leqslant j\leqslant g_1+g_2=g$, and $\omega^{(1)\prime}(x,p)$ and
$\omega^{(2)\prime}(x,p)$ are canonical differentials of second kind on $X_1$ and $X_2$ respectively.
\end{prop}
We refer the readers to \cite{fay} for standard terminology of canonical differentials of second kind and just remark for example
that $\omega^{(1)\prime}(x,p)$ is evaluated at $p$ with respect to the local coordinate $U_1$.

\ms

\ms
\ni{\bf 2.5} In the case of non-separating node $p$, $X_t$ degenerates into a stable curve $X_o$ with node at $p$, which  can be considered
as a connected Riemann surface $\hX_o$ with two points $a, b\in \hX_o$ identified.  Again, there exists small coordinate
neighborhoods $U_a, U_b$ of $a$ and $b$ respectively, and $\hX_o$ is the normalization of
$X_o$.  We may regard $U_a$ and $U_b$ as disks of fixed radius $\delta$ in some local coordinates around $a$ and $b$ respectively. The analytic structure can be given as in {\bf 2.4}. Let $0<\rho<1$.
Denote by
$\rho U_a$ and $\rho U_b$ disks of radius $\rho\delta$.

\begin{prop}\label{nonsep} (\cite{wentworth} page 437, \cite{fay} page 51, \cite{yamada} page 135)
We can find a normalized basis of $\Gamma(X_t,K_{X_t})$ for $t$ sufficiently close to $o$ such that for $x\in \hX_o-\rho U_a-\rho U_b,$
\begin{equation*}\begin{array}{cccl}
\omega'_i(x,t)&=&\omega'_i(x)-t[\omega'_i(b)\omega'(x,a)+\omega'_i(a)\omega'(x,b)]+\cO(t^2), & (1\leqslant i\leqslant g-1)\\
\omega'_g(x,t)&=&\omega'_{b-a}(x)-t[\gamma_1\omega'(x,b)+\gamma_2\omega'(x,a)]+\cO(t^2), &
\end{array}
\end{equation*}
where $\gamma_i$'s are some constants.  Moreover, the expression $\lim_{t\to 0}\cO(t^2)/t^2$ is a meromorphic form with poles only at $a$ and $b$, and
the coefficients has uniform convergence on $\hX_o-\rho U_a-\rho U_b$.
\end{prop}
In the above, $\omega'_{b-a}(z)=\frac1{2\pi i}\pd_z\log\frac{E(z,b)}{E(z,a)}$ and $E(z,a)$ is the prime form of $\hX_o$.   Since $E(z-a)$ in local
coordinates is given by $z-a$, we conclude that $\omega'_{b-a}(z)=\frac1{2\pi i}(\frac1{z-b}-\frac1{z-a})$ in local coordinates.

\ms
\ni{\bf 2.6} We recall the definition of generalized Weierstrass points on a projective algebraic curve as given in \cite{laksov} and \cite{og}.
A point $p\in X$ is called a Weierstrass point of the holomorphic line bundle $L$ (represented by $z$ as above) if there is an $s\in \Gamma(X,L)$ whose vanishing order at $p$ is at least $h^0(L):= \dim_{\mathbb{C}} \Gamma(X,L)$.
As in the case of the usual Weierstrass points, the Weierstrass points of a line bundle can also be defined in terms of the Wronskian of a basis of sections of
$L$ in \cite{laksov} and \cite{og}.

Consider the case that $X$ is a smooth curve of genus $g\geqslant 2$.  
Denote by $J_{d}$ the Picard variety of degree $d$.  Denote by 
$\Theta$ the theta divisor of $X$ in $J_{g-1}$.  There is a mapping $f_n:X\times \Theta\rightarrow J_{g-1+n}$ defined by
$f_n(x,\theta)=nx+\theta$.  
Then it is well-known that 

$x$ is a Weierstrass point of the line bundle $z$ if and only if 
\begin{equation}\label{theta}
    z=f_n(x,\theta)
\end{equation} for some 
$\theta\in \Theta$, which was taken as definition in \cite{neeman}.

Let $p$ be a Weierstrass point of $L$ over $X$. Then the weight of $p$, denoted $w_L(p)$, is defined as follows:
\vspace{1mm}
Let $s_1,...,s_m$ be a basis of $\Gamma(X,L)$ with distinct vanishing orders $\alpha_1<\cdot\cdot\cdot<\alpha_m$ at $p$, then 
$$w_L(p):= \alpha_1 +\cdot\cdot\cdot + \alpha_m -0-1-2-3-\cdot\cdot\cdot -(h^0(L)-1)$$
Notice that non-Weierstrass points have weight $0$.
Denote by $W(L)$ the set of all Weierstrass points of $L$ over $X$.

Let $h:X\rightarrow \bR$ be a continuous function on $X$. Let $L$ be any holomorphic line bundle of degree $g-1+m$ over $X$ ($m > g-1$).

Define the distribution
\begin{equation} \label{wm}
\mu^W_{X,L}:=\frac{\sum_{p\in W(L)} w_L(p)\cdot \delta_p}{\sum_{p\in W(L)}w_L(p)}
\end{equation}
where $\delta_p$ is the Dirac Delta at $p$. In case that there is no danger of confusion, we would simply denote
$\mu^W_{X,L}$ by $\mu^W_{L}$.
Then
$$
\int_X h \cdot \mu^W_L =\frac{\sum_{p\in W(L)} h(p)\cdot w_L(p)}{\sum_{p\in W(L)} w_L(p)}=\frac{1}{gm^2}\cdot\left(\sum_{p\in X} h(p)w_L(p)\right).
$$

\ms
\ni{\bf 2.7}  Recall the following result of \cite{neeman}, see also \cite{mumford1}. 
\begin{prop} 
Let $X$ be a Riemann surface of genus $g\geqslant 2$. Let $h:X\rightarrow \bR$ be a continuous function on $X$. Let $L$ be any line bundle of degree $g-1+m$ over $X$ ($m > g-1$). Then
$$ \frac{\sum_{p\in W(L)} h(p)\cdot w_L(p)}{\sum_{p\in W(L)} w_L(p)}=\int_X h \cdot \mu^W_{L} $$
converges to the constant
$$\frac{\int_X h \cdot \left(\omega_1\wedge\overline{\omega_1}+\cdot\cdot\cdot+\omega_g\wedge\overline{\omega_g}\right)}{\int_X \left(\omega_1\wedge\overline{\omega_1}+\cdot\cdot\cdot+\omega_g\wedge\overline{\omega_g}\right)}=\frac{1}{g}\cdot\int_X h\cdot \mu^B_X$$
as $m\to\infty$.\\\\
In the above, $\{\omega_1,...,\omega_g\}$ is an orthonormal basis of $\Gamma(X,K_X)$ and $\mu^B_X$ is the Bergman measure on $X$. 
\end{prop}

\ms
\ni{\bf 2.8} Let $P_{g,d}$ be the variety consisting of pairs $[C,L]$, where $C\in\cM_g$ and $L$ is a line bundle on $C$ of degree $d$.  In general, it is a subtle problem to have a natural canonical compactification of $P_{g,d}$ sitting above $\ocM_g$.  The difficulty is shown by the non-uniqueness of
extension of line bundle in the following example.  Consider a one-parameter family
of stable curves $\pi:\cC\rightarrow \Delta$ with $\Delta^*=\Delta-\{0\}\subset \cM_g$, where fibers 
$C_t, t\in \Delta$ is smooth and $C_{0}$ is nodal consisting of two components $C_{01}$ and $C_{02}$ meeting at a point $p$.  Let ${\cL}$ be a line bundle on $\pi^{-1}(\Delta^*)$ so that $L_t=\cL|_{C_t}$ is a line bundle on $C_t$ for $t\in\Delta^*$.  Then the extension of $\cL$ over $\Delta$ is not unique, since $\cL+\cO_{\cC}C_{01}$ would give another possible
extension apart from a given extension $\cL$ over $\Delta$.  Here $(\cL+\cO_{\cC}C_{01})|_{C_{02}}=(\cL+(p))|_{C_{02}}$ has degree $\deg(\cL|_{C_{02}})+1$
and $(\cL+\cO_{\cC}C_{01})|_{C_{01}}=(\cL-(p))|_{C_{02}}$ has degree $\deg(\cL|_{C_{02}})-1$.

For the case of irreducible stable $C_0$, the problem of compactification of $P_{g.d}$ is resolved 
by considering torsion-free coherent sheaves of rank one as given by \cite{ds}, and the above difficulty of uniqueness in extension 
does not occur since there is only one irreducible component.  In particular,
a line bundle with a fixed degree on $C_0$ extends to a line bundle on $C_t$ for
$t\in U$, a neighborhood of $0$ in $\ocM_g$. 

In the example above with two irreducible components meeting at a point, the problem of compactification was resolved in \cite{c}, in which the extension is
unique by considering 
line bundles of appropriate bidegree in $\Pic^{(d_1,d_2)}(C_0)=\Pic^{d_1}(C_{01})\times \Pic^{d_2}(C_{02})$, where the choice of bidegree is finite.  Recall that a nodal curve is {\it of compact type} if every node is separating.  In such a case, once we fix a multi-degree corresponding to a choice in {\cite{c}}, the extension is unique.  In particular, a line bundle with a fixed multi-degree on $C_0$ extends to a line bundle on $C_t$ for $t\in U$, a neighborhood of $0$ in $\ocM_g$.

In this article we study stable curves which are either irreducible or of compact type 
and consider line bundles which extend to a neighborhood $U$ of $0$ in $\ocM_g$.

\bs
\begin{center}
{\bf 3. Convergence of Bergman measure on a family of curves}
\end{center} 

\bs
\ni{\bf 3.1} {\bf Proof of Theorem 1}  

Let us first give a short outline of proof of Theorem 1a. 
It is well-known that in the setting of Theorem 1(a), a holomorphic one-form on
$X_t$ gives rise to a one form with at most a log pole at the node $p$.  Recall that the total residue of a meromorphic one-form on a connected Riemann surface is trivial.  If $p$ is separable so that $X_o$ consists of two irreducible components $X_1$ and $X_2$ of genus $a$ and $g-a$ respectively, the residue argument as above applied to the normalization $\hX_i$ of each component $X_i, i=1,2$, implies that a meromorphic one form cannot have
pole at a single point and hence the form is actually holomorphic.  In this case,  the sum of the Bergman kernels on $X_1$ and $X_2$ is precisely the limit
of the Bergman kernel on $X_t$.
If $p$ is non-separable, this corresponds to a meromorphic one form with a single pole at $p_1, p_2$ of opposite residues, where $\{p_1,p_2\}=\tau^{-1}(p)$.  In such case, there is a $g-1$ dimensional
space of holomorphic one-forms and one meromorphic one form with a log pole at
the node on $X_o$ coming from the convergence of the space of holomorphic one-forms from $X_t$.  One expects that the Bergman kernel of $X_t$ approaches the
Bergman kernel of $X_o$ as $t\rightarrow o$. It is however a bit tedious to
describe the convergence of the Bergman kernel since orthonormality 
is imposed in the definition of Bergman kernel as given in {\bf 2.2} and a log pole is not $L^2$-integrable.  We provide some details below.

\ms
\ni{\bf 3.2} {\it Theorem 1(a)(i)---}
This is already observed in Lemma 6.9 of [We].  For completeness of presentation, we explain the reason here parallel to our argument for (ii) in {\bf 3.3}.
We are given a family of curves $\pi:\cX\rightarrow S$ with $o\in S$ representing $X_o = X_1\cup X_2$. Let $\{\omega^{(1)\prime}_1(x,0),\dots,\omega^{(1)\prime}_{g_1}(x,0)\}$ and $\{\omega^{(2)\prime}_{g_1+1}(x,0),\dots,\omega^{(2)\prime}_{g_1+g_2}(x,0)\}$ be canonically normalized bases with respect to symplectic homology bases on $X_1$ and $X_2$ respectively.\\
On $X_1$, let $\{\omega^{(1)}_1(x,0),\dots,\omega^{(1)}_{g_1}(x,0)\}$ be an orthonormal basis of $\Gamma(X_1, K_{X_1})$ with respect to the natural $L^2$ norm as defined in \textbf{2.1}. Similarly for $\{\omega^{(2)}_{g_1+1}(x,0),\dots,\omega^{(2)}_{g_1+g_2}(x,0)\}$. Let $H_1(0), H_2(0)$ be the transformations so that 
$$ \omega^{(1)}_i(x,0) = H_1(0)_{ij} \cdot \omega^{(1)\prime}_j(x,0) \ \mbox{ for } \  1\leqslant i, j \leqslant g_1 $$
$$ \omega^{(2)}_j(x,0) = H_1(0)_{ij} \cdot \omega^{(1)\prime}_j(x,0) \ \mbox{ for } \  g_1+1 \leqslant i, j \leqslant g_1 + g_2$$

Let $H(t)$ be the transformation (for $|t|$ sufficiently small) such that 
\begin{equation} \label{0}
\begin{split}
\omega^{(1)}_i(x,t) = H(t)_{ij} \cdot \omega^{(1)\prime}_j(x,t) \ \mbox{ for } \  1\leqslant i, j \leqslant g_1 \\
\omega^{(2)}_i(x,t) = H(t)_{ij} \cdot \omega^{(2)\prime}_j(x,t) \ \mbox{ for } \  g_1+1 \leqslant i, j \leqslant g_1 + g_2
\end{split}
\end{equation}
and $H(t)^{-1}$ exists. Indeed,
\begin{equation}
H(t)=\left(\begin{array}{cc}
H_1(t) & 0\\
0 & H_2(t)
\end{array}\right)
\end{equation} where $H_1$ and $H_2$ are square matrices of size $g_1$ and $g_2$ respectively.
 Since $$\mu^B_{X_t}=\sum_{i=1}^g\omega_{X_t,i}\wedge\overline{\omega_{X_t,i}} = \sum_{i,j,k} H(t)_{ij} \overline{H(t)_{ik}} \omega'_j(x,t) \overline{\omega'_k(x,t)},$$ taking limit on both sides yields the result.

\ms
\ni{\bf 3.3} {\it Theorem 1(a)(ii)---} We consider degeneration of the Weierstrass points for a stable nodal curve with a node $p$ and degeneration
as given in Proposition \ref{nonsep}.  Hence we have a family of curves $\pi:\cX\rightarrow S$ with $o\in S$ representing $X_o$. 
$\hX_o$ is the normalization of $X_o$.  Let $\{\omega'_1(x,0),\dots,\omega'_{g-1}(x,0)\}$ be a canonically normalized basis with respect to
a symplectic homology basis  on $\hX_o$.

On $\hX_o$, we let $\{\omega_1(x,0),\dots,\omega_{g-1}(x,0)\}$ be an orthonormal basis of $\Gamma(\hX_o,K_{\hX_o})$ with respect to the 
natural $L^2$ norm as defined in \textbf{2.1}.  Let $J(0)$ be the transformation so that
\begin{equation}
\omega_i(x,0)=J(0)_{ij}\omega'_j(x,0) \ \mbox{ for } \  1\leqslant i, j\leqslant g-1.
\end{equation}
Let $\{\omega'_i(x,t)\}_{i=1}^g$ be the set of one forms on $X_t$ given by Proposition \ref{nonsep}.  There exists an invertible transformation $J(t)$ (for $|t|$ sufficiently small) satisfying
$$\omega_i(x,t)=J(t)_{ij}\omega'_j(x,t)  \ \mbox{ for } \  1\leqslant i, j\leqslant g-1$$ and $\{\omega_i\}_{i=1}^{g-1}$
being an orthonormal basis of 
$\span(\omega'_1,\dots,\omega'_{g-1})\subset \Gamma(X_t, K_{X_t})$.
Adding one more one form $\omega_g(x,t)$ so that $\{\omega_i\}_{i=1,\dots,g}$ gives an orthonormal basis of $\Gamma(X_t, K_{X_t})$, it follows that we can find a transformation $H$ containing $J$ as a submatrix so that 
\begin{equation} \label{1}
\omega_i(x,t)=H(t)_{ij}\omega'_j(x,t) \ \mbox{ for } \  1\leqslant i, j\leqslant g.
\end{equation}
Indeed, 
\begin{equation}
H(t)=\left(\begin{array}{cc}
J(t)&0\\
a^T(t)&b(t)
\end{array}\right)
\end{equation}
with $a^T(t)=(a_1(t),\cdots,a_{g-1}(t))$.  It follows that the inverse of $H$ is given by
\begin{equation} \label{2}
H^{-1}=\left(\begin{array}{cc}
J^{-1}(t)&0\\
-\frac{1}{b}a^T\cdot J^{-1}(t) &1/b(t)
\end{array}\right)
\end{equation}

From (\ref{1}), we know that
\begin{equation} \label{3}      
\omega_g(x,t)=\sum_{i=1}^{g-1}a_i(t)\omega'_i(x,t) + b(t)\omega'_g(x,t).
\end{equation}

By construction, $a_i(t)$ is smooth in $t$ for $1\leqslant i\leqslant g-1$.  Moreover, $\omega'_i(x,t)$ is uniformly bounded on $X_o$ when $t=0$ for $1\leqslant i\leqslant g-1$, and the expression
varies smoothly with respect to $t$.  Hence the expression $\sum_{i=1}^{g-1}a_i(t)\omega'_i(x,t)$ above is uniformly bounded for small $t$. It remains to estimate the term $b(t)\omega'_g(x,t)$.\\

Now, since $\omega_g \perp \span\{\omega'_1,...,\omega'_{g-1}\}$, we have
$$0=\int_{X_t}\omega_g(x,t)\wedge\overline{\sum_{i=1}^{g-1}a_i(t)\omega'_i(x,t)},$$
plugging (\ref{3}) into the above gives
\begin{equation} \label{4}
0=\int_{X_t}(\sum_{i=1}^{g-1}a_i(t)\omega'_i(x,t))\wedge\overline{\sum_{i=1}^{g-1}a_i(t)\omega'_i(x,t)}+\sum_{i=1}^{g-1}b(t)\overline{a_i(t)}\int_{X_t}\omega'_g(x,t)\wedge\overline{\omega'_{i}(x,t)}.
\end{equation}

We claim that for $1 \leqslant i \leqslant g-1$, there is the estimate
\begin{equation*}
\int_{X_t}\omega'_g(x,t)\wedge\overline{\omega'_i(x,t)}=o(t).
\end{equation*}
This follows from smoothness of $\pi$ and
\begin{equation}
\int_{X_o}\omega'_g(x,0)\wedge\overline{\omega'_i(x,0)}=0
\end{equation}
where $\omega'_g(x,0)=\omega'_{b-a}(x)$. The above identity is true because from our assumption, $\omega'_i(x,0)$ for $1\leqslant i\leqslant g-1$ is dual to a symplectic basis $\{A_i\}_{i=1,\dots,g-1}$. Hence $\int_{A_i}\omega'_g(x,0)=0$ 
for $1\leqslant i\leqslant g-1$ from the normalization in \textbf{2.2}
and so the claim is valid.

It follows from the claim and (\ref{4})
that 
\begin{equation}
\int_{X_t}(\sum_{i=1}^{g-1}a_i(t)\omega'_i(x,t))\wedge\overline{\sum_{i=1}^{g-1}a_i(t)\omega'_i(x,t)} = o(t)
\end{equation}

\ms
Recall that Proposition \ref{nonsep} gives rise to
\begin{equation}
\omega'_g(x,t)=\omega'_{b-a}(x)-t[\gamma_1\omega'(x,b)+\gamma_2\omega'(x,a)]+\cO(t^2),
\end{equation}
for $x\in \hX_o-\rho U_a-\rho U_b$ and the estimates in $t[\gamma_1\omega'(x,b)+\gamma_2\omega'(x,a)]$ and $\cO(t^2)$ are uniform.
Hence for fixed $\rho>0$, given any small $\epsilon>0$, we know that 
\begin{equation}
|\omega'_g(x,t)-\omega'_{b-a}(x)|<\epsilon
\end{equation}
if $t$ is sufficiently small and $0< t < \rho$.
Now 
\begin{eqnarray} \label{5}
\Vert \omega'_{b-a}\Vert_{\hX_o-\rho U_1-\rho U_2}^2&:=&\int_{\hX_o-\rho U_1-\rho U_2}\omega'_{b-a}\wedge\overline{\omega'_{b-a}}\nonumber \\
&\geqslant&\int_{(U_1-\rho U_1)\cup(U_2-\rho U_2)}\omega'_{b-a}\wedge\overline{\omega'_{b-a}}\nonumber \\
&\geqslant&c|\log\rho|
\end{eqnarray}
for some constant $c>0$ from direct integration.

Since $\Vert \omega_g(\cdot,t)\Vert_{X_t} =1$, and $\sum_{i=1}^{g-1}a_i(t)\omega'_i(x,t)$ is uniformly bounded for small $t$, it follows from
identity (\ref{3}) and the estimate (\ref{5}) that $b(t)>c_1|\log\rho|$ for some constant $c_1>0$ if $t<\rho$.
Hence
the one form $\omega_g(x,t)$ converges to $0$ on compacta  on $X_o-\{p\}\cong \hX_o-\{a,b\}$ as $t\rightarrow 0$.

Hence for $x\in \hX_o-\{p\}$, Proposition \ref{nonsep} implies that $\omega_{X_t,i}(x)\rightarrow \omega_{\hX_o,i}(x)$
for $1\leqslant i\leqslant g-1$ as $t \to 0$.  Since $\mu^B_{X_t}=\sum_{i=1}^g\omega_{X_t,i}\wedge\overline{\omega_{X_t,i}}$, 
it follows from the last paragraph that the limit of $\omega_{X_t,g}\wedge\overline{\omega_{X_t,g}}$ would concentrate
at the node $p$ as $t \to 0$. Here we note that the total measure 
$$\int_{X_t} \mu^B_{X_t}=\sum_{i=1}^g\Vert \omega_{X_t,i}\Vert^2=g$$
and $\sum_{i=1}^{g-1}\Vert \omega_{\hX_o,i}\Vert^2=g-1.$  The discrepancy is precisely given by the delta measure at
the point $p$, since the mass cannot be concentrated anywhere else according to the discussions above. 

\ms
\ni{\bf 3.4}
{\it Theorem 1(b)---}  
This follows from {\bf 3.2, 3.3} and induction.  Suppose $k=2$.  Suppose $X_o$ is a stable curve with two nodes obtained after contracting two nodes from families of smooth curves, corresponding to a point at the
boundary of the Deligne-Mumford compactification $D_i\cap D_j\subset D=\ocM_g-\cM_g$ for some $i\neq j$.  We may consider a
local two dimensional holomorphic family of curves 
$X(s,t)$ for 
$(s,t)\in \Delta\times \Delta$, so that $X(s,t)$ is smooth for $s\neq0$ and $t\neq0$, $X(0,t)\in D_i$ and $X(s,0)\in D_j$, and $X(0,0)=X_o$.  

We consider
first a point $X(0,t)$ at $\pd\cM_g$ obtained by
contracting $1$ real cycle giving rise to a node $p_1(t),$ which may be assumed to be a fixed node $p_1$ with respect to a local 
trivialization of the family.  This is obtained by letting
$s\rightarrow 0$ in $X(s,t)$.  Let  $\hX(0,t)$ be the normalization of $X(0,t)$.  Theorem \ref{bergmanlimit} implies that $$\lim_{s\rightarrow0}\mu^B_{X(s,t)}=\mu^B_{X(0,t)}=\mu^B_{\hX(0,t)}+\delta_{p_1}.$$
$X_o$ is obtained by contracting a real $1$-cycle on $X(0,t)$ to a node $p_2$, which corresponds to contracting a real cycle on $\hX(0,t)$ to 
a node $\hp_2$, where $\hp_2$ corresponds exactly to
the node $p_2$ on $X_o$.  Now the normalization of $\hX(0,0)$ is precisely $\hX_o$.  Hence Theorem \ref{bergmanlimit}a again implies that
$$\lim_{t\rightarrow0}\mu^B_{\hX(0,t)}=\mu^B_{\hX(0,0)}+\delta_{p_2}=\mu^B_{\hX_o}+\delta_{p_2}.$$

Combining the above two identities, we see that 
$$\lim_{t\rightarrow0}\lim_{s\rightarrow0}\mu^B_{X(s,t)}=\mu^B_{\hX_o}+\delta_{p_2}+\delta_{p_1}.$$

Note that the arguments of \cite{fay}, \cite{yamada} and \cite{wentworth} concerning behavior of period matrices corresponding to contraction of a real $1$-cycle applies equally well to a family of degenerating curves obtained by contracting two different non-intersecting real $1$ cycles as well.  The end result depends only on $X_o$ and is independent of the paths of degeneration taken.  

Hence Theorem 1(b)
is proved for $k=2$.
The same proof clearly works for $k>2$ as well. 

\qed

\bs
\begin{center}
{\bf 4. Convergence of Weierstrass measure on a family of curves} 
\end{center} 



\ms
\ni{\bf 4.1}
Suppose $C$ is a stable curve with nodal singularities at $z_i$, $i=1,\dots,n$.  Let $\pi:\tC\rightarrow C$ be the
normalization of $C$ so that $\pi^{-1}(z_i)=\{a_i,b_i\}$.
Let $U$ be a small coordinate neighborhood of $z_i$.  Then the dualising sheaf $\omega_C$ is generated by holomorphic $1$-forms in
a neighborhood of a regular point on $C$ or $\tC$, and by meromorphic $1$-forms $\eta$ with at worst simple poles at $a_i, b_i$ over $U$, satisfying
$$\Res_{a_i}(\eta)+\Res_{b_i}(\eta)=0.$$
Now let $L$ be an ample line bundle on a stable curve $C$.
Let $\psi, \tau$ be the generators of $L$ and $\omega_C$ over $U$ respectively.  Let $n=h^0(C,L)$ and $\phi_1,\dots,\phi_n$ be a basis of $H^0(C,L)$.  Define
$F_{i,j}\in\Gamma(U,\cO_C)$ inductively by
\begin{eqnarray*}
F_{1,j}\psi &:=&\phi_j|_U \hspace{6mm} j=1,\dots,n,\\
F_{i,j}\tau &:=&dF_{i-1,j} \hspace{6mm} i=2,\dots,n,\ j=1,\dots,n
\end{eqnarray*}
Define also
$$\rho=\det(F_{i,j})\psi^n\tau^{(n-1)n/2}.$$
It follows easily by checking compatibility on different
charts that $\rho$ defines a section in $H^0(C,L^n\otimes \omega_C^{(n-1)n/2})$.  Then as in \cite{widland-lax}, we define
$p\in C$ to be a Weierstrass point of $L^n\otimes \omega_C^{(n-1)n/2}$ if and only if 
\begin{equation}\label{weierstrass}
    \ord_p\rho>0. 
\end{equation}

It follows from \cite{widland-lax} that the number of Weierstrass points counted with multiplicity is given by
$n\deg(L)+(n-1)n(g-1)$. Suppose $n>g-1$, the Riemann-Roch formula shows that $\deg(L)=g-1+n$ and hence the number of 
Weierstrass points counted with multiplicity is then given
by $n^2g$.

\ms
\ni{\bf 4.2} Let us consider first in details the situation that $\hX_o$ has
genus $0$, a case that partly motivates the present paper.

\begin{lemm}
Let $X_o$ be an irreducible rational nodal curve.  For each $m\in \bN$, $\lim_{t\rightarrow0}\mu^W_{X_t,mK_t}=\mu^W_{X_o,mK_o}$.
\end{lemm}

\ni{\bf Proof} We assume that $X_o$ is a rational curve with $g$ double points. Hence $X_o$ is formed by identifying $g$ pairs of distinct points $b_i$ and $c_i$, $i=1,\dots,g$ on $P_{\bC}^1$.
In this case, the discussions in {\bf 5.1} could be realized concretely as
follows (for details, see \cite{mumford2},\cite{lax}).\\\\
The dualizing sheaf of $X_o$ is spanned by
$$ \omega_i = \frac{dz}{z-b_i} - \frac{dz}{z-c
_i} \hspace{6mm} i=1,2,...,g$$
Thus the period lattice $\Lambda$ is generated by the $g$ vectors
$$\{(2\pi\sqrt{-1},0,...,0),(0,2\pi\sqrt{-1},0,...,0),...,(0,...,0,2\pi\sqrt{-1})\}$$
Hence the generalized Jacobian is $\mathbb{C}^g/\Lambda \cong (\mathbb{C}^*)^g$.
Let $X_o^s$ be the set of points on $X_o$ with all $b_i, c_i$ removed, i.e., the set of smooth points. And we further assume $x_o = \infty\in X_o^s$, this can be done by choosing appropriate coordinate.
Choosing $x_o$ as basepoint, we define the Abel mapping $\varphi:X_o^s \to J(X_o)\cong (\mathbb{C}^*)^g$ by
$$ \varphi(x) = \left( exp \left(\int_{x_o}^x \omega_1\right),\cdots, exp \left(\int_{x_o}^x \omega_g  \right) \right) $$
$$ \hspace{-15mm} = \left( \frac{x-b_1}{x-c_1},\cdots, \frac{x-b_g}{x-c_g} \right) $$
This induces a map $\varphi:(X_o^s)^{m}\rightarrow (\bC^*)^g$
given by 
$$\varphi(\sum_k n_kx_k)=\left( \prod_k\left(\frac{x_k-b_1}{x_k-c_1}\right)^{n_k},\cdots,\prod_k \left(\frac{x_k-b_g}{x_k-c_g}\right)^{n_k}\right)$$
Let $\lambda_i = exp(z_i)$ be coordinates on $(\mathbb{C}^*)^g$. 
Define $\tau$ on $J$ by 
$$\tau(\lambda_1,\dots,\lambda_g)=\det
\begin{bmatrix}
1-\lambda_1&\cdots&1-\lambda_g\\
b_1-c_1\lambda_1&\cdots&b_g-c_g\lambda_g\\
\vdots&&\vdots\\
b_1^{g-1}-c_1^{g-1}\lambda_1&\cdots&b_g^{g-1}-c_g^{g-1}\lambda_g
\end{bmatrix}$$
Let $X_t$ be a family of smooth curves of genus $g$ degenerating to $X_o$. It is a standard fact (\cite{fay}) that the period matrices $\Omega_{ij}(t)$ of $X_t$ satisfy
$$ Im(\Omega_{ii}(t)) \to \infty \hspace{8mm} \text{ as } t\to 0 $$
and 
$$ \Omega_{ij}(t) \text{ are continuous } \hspace{4mm} \text{ for } |t|<\epsilon $$
Let $\Omega_{ii}(t)$ be the diagonal of $\Omega(t)$, we have, upon direct computation, that, as $t\to 0$, 
$$\theta\left(z-\frac{1}{2}\Omega_{ii}(t),\Omega(t)\right)$$ converges to $\tau$ up to a constant multiple. (\cite{mumford2}, page 3.253.)\\\\
Hence $\tau$ defines the Jacobian divisor $\theta_o$ on $X_o$.\\\\
To sum up, consider a family of stable curves $X_t$ with smooth $X_t$ when $t\neq0$ and $X_o$ is a stable curve with double points in the sense of Deligne-Mumford.  In the case that $X_o$ is just a rational curve with double points, we know that there is a convergence of $\theta_t$ to $\theta_o$. Hence from the definition of Weierstrass points earlier, there is a convergence of the Weierstrass divisors as claimed in the statement of the lemma.

\qed

\ms
\ni{\bf 4.3} In this subsection, we generalize the argument in the previous subsection
to the case of arbitrary genus.

Consider now the family of stable algebraic curves of
genus $g$, $\pi:\cX\rightarrow S$ as before so that fibers $X_t$ are smooth except for $X_o$ which has nodal singularity with a single node
at $p$.  We define $\mu^W_{mL_t}$ as the distribution associated to Weierstrass points 
on $X_t$ as in (\ref{wm}).

\begin{lemm}  Assume that $X_o$ is stable.   Then
for each $m\in \bN$,\\ $$\lim_{t\rightarrow0}\mu^W_{X_t,mL_t}=\mu^W_{X_o,mL_o}.$$
\end{lemm}

\ms\ni{\bf Proof} Since there is only one node, $X_o$ is either irreducible or of compact type.
 A Weierstrass point on a stable curve is defined by (\ref{weierstrass}).
For the case of irreducible stable curve, the lemma follows from \cite{lax2} Theorem 1.   In the case that $X_t$ is stable and of compact type, it follows from
\cite{en} Theorem 8.4, \cite{es} Theorem 6.  In either case, this follows from
the convergence of the Wronskian in the definition of Weierstrass points.


An alternative approach closer to the description in (4.2) for irreducible stable rational curve can be given as follows, making use of the alternative
definition of Weierstrass points as given in (\ref{theta}). Denote by
$\Pic^d=J^d$ the Picard variety of degree $d$ on $X$, parametrizing line bundles of degree $d$.  $\Pic^d$ may not be projective if $X$ is not smooth.  In such case, we may consider compactified Jacobian and theta divisor as defined in \cite{alexeev1}, \cite{alexeev2}.  
If $X$ is a smooth curve, the theta divisor can be defined intrinsically as the locus of $L\in \Pic^{g-1}(X)$ with $h^0(X,L)\neq0$.  This is used as definition for stable curve as well in the following way.

According to \cite{alexeev1}, there is a complete moduli of semiabelic pairs defined in \cite{alexeev1}, for which the compactified Jacobian and its theta divisor is such a pair. From the work of Simpson in 1.21 of \cite{simpson}, see also the explanation in 1.2-1.3 of \cite{alexeev2}, there is a family $\pi:\mathcal{J}\rightarrow S$ of compactified Jacobians of degree $g-1$ over the base curve $S$ in which each fiber is the compactified Jacobian of $X_t$, for which it follows from \cite{alexeev2} that we
may choose an arbitrary polarization.
Let $\mathcal{L}$ be the ample line bundle which gives the polarization.  From \cite{alexeev2}, $\pi_*\mathcal{L}$ is invertible by cohomology and base change. Choosing a trivialization of $\pi_*\mathcal{L}$, this gives a section $s\in H^0(S,\pi_*\mathcal{L})$ whose restriction $s_t$ to the fiber over $t$ is the unique section of $L_t$. This gives a family of theta divisors $\Theta_t, \forall t\in S$.

The setting above implies that the theta divisor on $J_t$ for $t\neq o$ converges to
$J_o$ as $t\rightarrow o$, which is a restatement of Mumford in higher genus case at the central fiber.  

Then from the definition of Weierstrass points
and weights, it follows that 
$$\lim_{t\rightarrow0}\mu^W_{X_t,mL_t}=\mu^W_{X_o,mL_o}.$$
\qed

\ms
We remark that the assumption that $X_t$ is stable irreducible or of compact type is
used in the second approach above.  It is well known that there is isomorphism between
Picard varieties of different degrees for smooth curves, after translation by a based line
bundle with degree the difference of the two.  This could also be done for stable irreducible curves or curves of compact type.  The problem is in general subtle for
arbitrary stable curves.  We refer the readers to \cite{gz} for some results in this
direction.

\bs
\begin{center}
{\bf 5. Weierstrass measure on stable curves} 
\end{center}

\ms
\ni{\bf 5.1} Recall that $p$ is the node on $X_o$, corresponding to $a$ and $b$ on $\hX_o$.
Let $U$ be a small neighborhood of $p$ corresponding to the union of two disks $U_a$ and $U_b$ around $a$ and $b$ respectively.
Let $\Delta$ be a sufficiently small neighborhood of $o$ in $S$.  We may assume that $U$ can be extended to $\cU$ on $\pi^{-1}(\Delta)\cap \cX$
so that for $U_t=\cU\cap X_t$,  $X_t-U_t$ is diffeomorphic to $X_o-U_o$ for all $t\in \Delta-\{o\}$.  The following lemma follows immediately from
the steps of proof in \cite{neeman}.

\begin{lemm}
\begin{equation}
\lim_{m\rightarrow\infty}\lim_{t\rightarrow o}\mu^W_{mL_t}(x)=\lim_{t\rightarrow o}\lim_{m\rightarrow\infty}\mu^W_{mL_t}(x)
\end{equation}
uniformly for $x\in (\pi^{-1}(\Delta)-\cU)$.
\end{lemm}

\ms
\ni{\bf Proof} We are going to follow the proof of Neeman in Chapter 2 of \cite{neeman} on $x\in (\pi^{-1}(\Delta)-\cU)$.  The reason that the argument
goes through is that we have nice convergence of the Bergman metric and geometry
on $(\pi^{-1}(\Delta)-\cU)$ as $t\rightarrow o$.  

As in \cite{neeman}, we consider $f_{t,n}:X_t\times \Theta_t\rightarrow J_{t,g-1+n}$ given
by $f_{t,n}(x,\theta)=nx+\theta$.
Let $F_t\subset X_t\times \Theta_t$ be the union of the singular set of 
$X_t\times\Theta_t$ and the ramification locus of $f_{t,n}.$  
By abuse of language, we denote by $J_{t,g-1+n}\cap((\pi^{-1}(\Delta)-\cU))$ the
subset of $J_{t,g-1+n}$ given by the image of the Jacobian image of $X_t\cap ((\pi^{-1}(\Delta)-\cU))$.
We are actually considering
only the restriction of $f_{t,n}$ to $f^{-1}(J_{t,g-1+n}\cap(\pi^{-1}(\Delta)-\cU))$, namely
$$f_{t,n}:X_t\times \Theta_t|_{f^{-1}(J_{t,g-1+n}\cap((\pi^{-1}(\Delta)-\cU))}\rightarrow J_{t,g-1+n}\cap((\pi^{-1}(\Delta)-\cU)).$$
Now for a translational invariant vector field $V_t$ on $J_{t,g-1+n}\cap((\pi^{-1}(\Delta)-\cU))$, we have
$f_{t,n}^{-1}(V_t)=V_{t,1}^n\oplus V_{t,2}^n$.  As in Lemma 2.3 of \cite{neeman}, we have
\begin{equation}
V_{t,1}^n=\frac1n V_{t,1}^1=\frac1n V_{t,1}, \  \ V_{t,2}^n=V_{t,2}^1=V_{t,2}.
\end{equation}
As in Lemma 2.3 of \cite{neeman}, the operators 
$f_{t,n}^{-1}(V_t^1)\cdots f_{t,n}^{-1}(V_t^1)$ on compact space $ D_t\subset(X_t\times \Theta_t-F_t)|_{f^{-1}(J_{t,g-1+n}\cap((\pi^{-1}(\Delta)-\cU))}$ are uniformly bounded as operators
$C_{D_t}^{r+m}\rightarrow C_{D_t}^r$ as $n$ varies and is uniform in $t\in \Delta$.
Let $h_t:X_t\times \Theta_t|_{f^{-1}(J_{t,g-1+n}\cap((\pi^{-1}(\Delta)-\cU))}\rightarrow \bR$ be a smooth function with compact support. Define $\Av^nh_t:J_{t,g-1+n}\rightarrow \bR$ by
$$(\Av^nh_t)(z)=\frac1{gn^2}\sum_{x\in W_t(z)}h_t(x,z),$$
where $W_t(z)$ denotes the set of Weierstrass points of $z$ with multiplicities.
Again, the argument of Lemma 2.6 of \cite{neeman} implies that there exists $M\in \bR$
such that for all $n>g-1$ and $t\in \Delta$ that 
$$\Vert V_t^1\cdots V_t^m(Av^n h_t)\Vert_\infty\leqslant M.$$
Similarly, it follows as Lemma 2.7 of \cite{neeman} that the
$k$-th Fourier coefficient of 
$\Av^nh_t$ for $k=(k_1,\dots,k_{2g})\neq (0,0,...,0)$, denoted by $(\Av^nh_t)\hat (k)$, satisfies
$$(\Av^nh_t)\hat (k)\leqslant\frac M{\sum_{i=1}^{2g}|k_i|^{2g+1}}.$$
Then Lemma 2.8 of \cite{neeman} implies that the $0$-th Fourier coefficient of $\Av^nh_t$, denoted by $(\Av^nh_t)\hat (0)$, satisfies
$$(\Av^nh_t)\hat (0)=\int_{C_t\times \Theta_t|_{f^{-1}(J_{t,g-1+n}\cap((\pi^{-1}(\Delta)-\cU))}} h_t \cdot dB_t.$$
The argument of Lemma 2.9 of \cite{neeman} then implies that
$$\Av^nh_t-\int_{C_t\times \Theta_t|_{f^{-1}(J_{t,g-1+n}\cap((\pi^{-1}(\Delta)-\cU))}}h_t \cdot dB_t$$ converges
uniformly to $0$ as $n\rightarrow\infty$ and uniformly in $t\rightarrow o$.
Now we may use the argument of Lemma 2.9 of \cite{neeman} to show that the above
convergence actually holds for arbitrary continuous function $h_t$ on 
$C_t\times \Theta_t|_{f^{-1}(J_{t,g-1+n}\cap((\pi^{-1}(\Delta)-\cU))}$, uniformly in $n\rightarrow \infty$ and
in $t\rightarrow o$.  In particular, we may interchange the order
of limits as given in our statement.

\qed
\medskip


\ms
\ni{\bf 5.2} 
{\bf Proof of Theorem 2}
Let $V$ be a small neighborhood of a node $p\in X_o$ as mentioned at the beginning of {\bf 5.1}.  We extend $V$ smoothly
to a neighborhood $\cV$ in the total family and use the same notation to denote
$\cV\cap X_t$ for $t$ sufficiently small.  Note that all of this can be performed
in a local coordinate as discussed in {\bf 2.4, 2.5}.
It follows from Lemma 2, Lemma 3 and Theorem \ref{bergmanlimit} that
for any small neighborhood $U$ of the node and any
$x\in U$, 
\begin{eqnarray}\label{id}
\lim_{m \to \infty}\mu_{mL_o} &=& \lim_{m\to\infty}\lim_{t\to o} \mu^W_{mL_t}\\
&=& \lim_{t\to o}\lim_{m\to\infty} \mu^W_{mL_t}\nonumber\\
&=& \frac1g \cdot \lim_{t \to o} \mu^B_{X_t}.\nonumber
\end{eqnarray}
Hence by shrinking $V$ to the node, we conclude that

\begin{equation} \label{18}
    \lim_{m \to \infty}\mu_{mL_o}|_{X_o-\{p\}} =1/g \cdot \lim_{t \to o} \mu^B_{X_t}|_{X_o-\{p\}}=1/g \cdot \mu^B_{X_o}|_{X_o-\{p\}}.
\end{equation}

Consider first the case (a) that the nodal point $p$ is separating.  In this case,
the sum of genera of the two components of $X_o$ is precisely $g$.
Hence 
\begin{equation}  \label{19}
\frac1g\int_{X_o}\mu^B_{X_o}=\frac1g\int_{X_o-\{p\}}\mu^B_{X_o}=1.\end{equation}

On the other hand,
\begin{equation}\lim_{m \to \infty}\mu_{mL_o}|_{X_o-\{p\}} =\lim_{m \to \infty}\lim_{t\rightarrow o}\mu_{mL_t}|_{X_o-\{p\}}
\end{equation}
and each
$\lim_{m \to \infty}\mu_{mL_t}$ is given by $\frac{1}{g}\cdot \mu^B_{X_t}$ and hence 
\begin{equation} \label{21}
\int_{X_o}\lim_{m \to \infty}\mu_{mL_o}=1.
\end{equation}
It follows from equations (\ref{18}), (\ref{19}) and (\ref{21}) that there is no mass trapped in $p$ and
hence 
\begin{equation}
    \lim_{m \to \infty}\mu_{mL_o}|_{X_o} =1/g \cdot \lim_{t \to o} \mu^B_{X_t}|_{X_o}=1/g\cdot \tau_*\mu^B_{\hat X_o}.
\end{equation}
This corresponds to (a).

Consider now (b) for which  $X_o$ is obtained from contracting a real
cycle on $X_t$ to a non-separating node on $X_o$.  The Genus of $X_t$ is $g$ for each $t\neq 0$ and the genus of $\hX_o$ is $g-1$.  It follows that  
\begin{eqnarray*}
\lim_{m \to \infty}\int_{X_t}\mu_{mL_t}&=&g \hspace{7mm} \text{for} \hspace{2mm} t\neq 0;\\
\lim_{m \to \infty}\int_{X_o}\mu_{mL_o}&=&g-1.
\end{eqnarray*}

It follows that the difference in the measure is supported at the node as
$t\rightarrow 0$.  Hence 
$$\lim_{m \to \infty}\mu_{mL_o}|_{X_o} =1/g \cdot (\lim_{t \to o} \mu^B_{X_t}|_{X_o}+\delta_p)=1/g\cdot (\tau_*\mu^B_{\hat X_o}+\delta_p).$$

\qed

\ms
\ni{\bf 5.3} 
{\bf Proof of Theorem 3}

 Let $X_o$ be a stable curve represented by a point on $\ocM_g-\cM_g$.  We may regard $X_o=X_0$ as the degeneration of a family
of smooth curves $X_t, t\in \Delta^*$ by contracting a finite number of
real $1$-cycles $\gamma_i, i=1,\dots,l$.  Denote by $p_i$ the nodes on $X_o$.

Consider first the case (b) that $X_o$ is of compact type.  In such case,
the normalization $\nu:\hX_o=\cup_{i=1}^{l+1}Y_i\rightarrow X_o$ has $l+1$ disconnected components
and there exist $q_{12}\in Y_1$, $q_{i1}, q_{i,2}\in Y_i$ for $2\leqslant i\leqslant l-1$, and $q_{l,1}\in Y_l$ so that $\nu(q_{i+1,1})=\nu(q_{i,2})=p_i$.
It is known that the sum of the genera satisfies $\sum_{i=1}^lg(Y_i)=g$.
Let $V$ be a
neighborhood of the nodes.  $V$ consists of several components if $l>1$.  In such case,
identity (\ref{id}) still holds and for $x\in X_o-V$ in the notation of proof of
Theorem 2,
\begin{equation}\label{ie}
\lim_{m\rightarrow\infty}\mu^W_{mL_o}=\frac1g \cdot \lim_{t \to o} \mu^B_{X_t}.
\end{equation}
From Theorem 1, we conclude that 
$\lim_{t \to o} \mu^B_{X_t}=\mu^B(\hX_0)=\sum_{i=1}^l\mu^B(Y_i)$, since all the nodes are
separating.  Hence we conclude that
\begin{equation}\label{if}\lim_{m\rightarrow\infty}\mu^W_{mL_o}=\frac1g(\sum_{i=1}^l\mu^B(Y_i))
\end{equation}
on $X_o-V.$
After shrinking $V$ to the
points $p_i$, we conclude that the identity (\ref{ie}) holds everywhere on $X_o-\{p_1,\dots,p_l\}$.  Since $\sum_{i=1}^lg(Y_i)=g$, we know that
$$\int_{X_o-\cup_{i=1}^l \{p_i\} }\mu^W_{mL_o}=\frac1g\int_{X_o-\cup_{i=1}^l \{p_i\} }\sum_{i=1}^l\mu^B(Y_i)=1,$$
where we used the fact that $Y_i$ are smooth for $1\leqslant i\leqslant l$ and hence $\mu^B(Y_i)$ are smooth measures as well. It follows that no mass is trapped in
$p_i, i=1,\dots,l$.  Hence the identity (\ref{if}) holds as a measure everywhere on $X_o$.

\ms
Consider now the case (a).  In this case, $X_o$ is of irreducible with 
$k$ non-separable nodes $p_i, i=1,\dots,k.$  The normalization
$\nu:\hX_o\rightarrow X_o$ is smooth and irreducible.  There are points $q_{ij}, j=1,2, i=1,\dots,k$ on $\hX_o$ such that $\nu(q_{ij})=p_i$.  The identity (\ref{ie}) still holds in this case.  Hence as in (b),  

$$g\lim_{m\rightarrow\infty}\mu^W_{mL_o}=\nu_*\mu^B_{\hX_o}+\mu_o,$$
where $\mu_o$ is supported on the nodes $\{p_1,\dots,p_l\}$.  
From our normalization, $\int_{X_o}\mu^W_{mL_o}=1$ and $\int_{X_o}\mu^B_{\hX_o}=g-l$, we conclude that
$\int_{X_o}\mu_o=l.$  Hence we may assume that $\mu_o=\sum_{i=1}^la_i\delta_{p_i}$ with $0\leqslant a_i$ and $\sum_{i=1}^la_i=l$.
In other words,
\begin{equation}\label{comp}
  g\lim_{m\rightarrow\infty}\mu^W_{mL_o}=\nu_*\mu^B_{\hX_o}+ \sum_{i=1}^la_i\delta_{p_i}. 
\end{equation}

We claim that $a_i\leqslant 1$ for $1\leqslant i\leqslant l$.
For simplicity of explanation, we consider first the case that 
$X_o$ has exactly two nodal points $p_1,p_2\in X_o$ in $\ocM_g$.  In our setting,
the line bundle $L_o$ at $X_o$ extends to a neighborhood $U$ of $o$ in $\ocM_g$.
Consider a deformation family of stable curves $X_t$ centered at $X_o$  so that each $X_t$
has precisely one nodal point $p_{1t}$ for $t\neq 0$ and $X_o$ corresponds to $t=0$.  In other words,
$p_2$ on $X_o$ is the result of contracting a real cycle $C_2$ on $X_t$. 
Here we may assume that $t\in \Delta$, a small disk centered at $o$ and that
coordinates described as in {\bf 4.2, 4.3} are used.

Since $\hX_t$ has just a single node at $p_{t1}$, we know that
$g\mu^W_{X_t,L_t}=g\lim_{m\rightarrow\infty}\mu^W_{X_t,mL_t}=(\nu_{X_t,p_{t1}})_*\mu^B_{X_t}+a_{t1}\delta_{p_{t1}}$ with $a_{t1}=1$ from Theorem 2, where
 $\nu_{\hX_t,p_{t1}}:\hX_t\rightarrow X_t$ is the normalization $X_t$ at $p_{1t}$ on $X_t$.  
 This holds for all $t\in \Delta^*$.  In particular, in taking $t\rightarrow0$, and applying Fatou's Lemma, we conclude that 
 $$a_1=a_{01}\leqslant \liminf_{t\rightarrow0}a_{t1}\leqslant 1.$$
 Similarly, $a_2\leqslant 1$.
 
 In the general situation of $k>2,$  we choose a local family of irreducible stable curves $X_t$ centered at $X_o$ such
 that $X_t$ has only a node at $p_{t1}$ and is smooth elsewhere for $t\in \Delta^*.$ The same argument as above shows that $a_1\leqslant 1.$  Applying the same argument to
 $p_{ti}$ for $1\leqslant i\leqslant k$, we conclude that $a_i\leqslant 1$ for all $i$ and hence
 the claim is proved.
 
 Since $\sum_{i=1}^la_i=l$, it follows from the claim that actually
 $a_i=1$ for $1\leqslant i\leqslant k$.
 
 \qed
 

 


\ms
\ni{\bf 5.4} {\bf Proof of Corollary 1}  

Let $X$ be a Riemann surface of genus $g$.  After contracting a non-separating real one cycle which is homologically non-trivial,
we obtain a stable curve $X_1$ with a node.  $X_1$ lies in the boundary of the Deligne-Mumford compactification $\ocM_g-\cM_g$.
The normalization of $\hX_1$ is a curve of genus $g-1$.  Repeat the
above procedure by contracting a non-separating real $1$-cycle on $\hX_1$.  It corresponds to contracting another non-separating
real $1$-cycle on $\hX_1$ and we arrive at a stable curve $X_2$ with
two separate nodes.  Inductively after $g$ steps, we arrive at a
rational curve $X_g$ with $g$ nodes.

Suppose now $X_o$ is a rational curves with nodal points obtained as 
above. Application of Theorem \ref{mainthm} to $X_o$ gives precisely the formula in Corollary 1.

\qed

\noindent {\bf Note Added in Proof } We would like to mention the recent work of \cite{am}, \cite{dj}, \cite{shivaprasad} on related topics, in particular, on tropical curves, made known to us after acceptance of the paper.



\begin{thebibliography}{99}

\bibitem[A1]{alexeev1}
V. Alexeev, 
Complete moduli in the presence of semiabelian group action. Ann. of Math. (2) 155 (2002), no. 3, 611--708, MR1923963, Zbl 1052.14017.

\bibitem[A2]{alexeev2}
V. Alexeev, 
Compactified Jacobians and Torelli map. Publ. Res. Inst. Math. Sci. 40 (2004), no. 4, 1241--1265, MR2105707, Zbl 1079.14019.

\bibitem[AK]{ak}
A. Altman, S. Kleiman, 
Compactifying the Picard scheme. Adv. in Math. 35 (1980), no. 1, 50--112, MR0555258, Zbl 0427.14015.
 
\bibitem[Am]{am} 
O. Amini, Equidistribution of Weierstrass points on curves over non-Archimedean fields, arXiv 1412.0926.
 
\bibitem[BG]{ballico-gatto}
E. Ballico, L. Gatto,
Weierstrass points on singular curves. Rend. Sem. Mat. Univ. Politec. Torino 55 (1997), no. 2, 145--170 (1998), MR1680491, Zbl 0923.14020

\bibitem[B]{burnol}
J.-F. Burnol, 
Weierstrass points on arithmetic surfaces. Invent. Math. 107 (1992), no. 2, 421--432, MR1144430, Zbl 0723.14019

\bibitem[C]{c} 
L. Caporaso, 
A compactification of the universal Picard variety over the moduli space of stable curves. J. Amer. Math. Soc. 7 (1994), no. 3, 589--660, MR1254134, Zbl 0827.14014

\bibitem[CE]{ce} 
L. Caporaso, E. Esteves, 
On Abel maps of stable curves. Michigan Math. J. 55 (2007), no. 3, 575--607, MR2372617, Zbl 1134.14019

\bibitem[CP]{cp} 
J. Coelho, M. Pacini, 
Abel maps for curves of compact type. J. Pure Appl. Algebra 214 (2010), no. 8, 1319--1333, MR2593665, Zbl 1194.14039

\bibitem[dJ]{dj} R. de Jong,  Faltings delta-invariant and semistable degeneration. J. Differential Geom. 111 (2019), no. 2, 241-301.

\bibitem[D]{delcentina} 
A., Del Centina, 
Weierstrass points and their impact in the study of algebraic curves: a historical account from the "Lückensatz'' to the 1970s. Ann. Univ. Ferrara Sez. VII Sci. Mat. 54 (2008), no. 1, 37--59, MR2403373, Zbl 1179.14001

\bibitem[DS]{ds} 
C. D'Souza, 
Compactification of generalised Jacobians. Proc. Indian Acad. Sci. Sect. A Math. Sci. 88 (1979), no. 5, 419--457, MR0569548

\bibitem[EN]{en}
Esteves, E., Nogueira, P.,
Generalized linear systems on curves and their Weierstrass points. Comm. Algebra 41 (2013), no. 3, 989--1016, MR3037176, Zbl 1271.14041

\bibitem[ES]{es}
Esteves, E., Salehyan, P., 
Limit Weierstrass points on nodal reducible curves. Trans. Amer. Math. Soc. 359 (2007), no. 10, 5035--5056, MR2320659, Zbl 1121.14016

\bibitem[F]{fay}
J. Fay,
Theta functions on Riemann surfaces. Lecture Notes in Mathematics, Vol. 352. Springer-Verlag, Berlin-New York, 1973. {\rm iv}+137 pp, MR0335789, 

\bibitem[FL]{little-furio}
K. A. Furio, J. B. Little,
On the distribution of Weierstrass points on irreducible rational nodal curves. Pacific J. Math. 144 (1990), no. 1, 131--136, MR1056669

\bibitem[GL]{garcia-lax}
A. Garcia, R. F. Lax,
Rational nodal curves with no smooth Weierstrass points. Proc. Amer. Math. Soc. 124 (1996), no. 2, 407--413, MR1322924

\bibitem[GZ]{gz}
S. Grushevsky, D. Zakharov,
The double ramification cycle and the theta divisor. Proc. Amer. Math. Soc. 142 (2014), no. 12, 4053--4064, MR3266977

\bibitem[HM]{hm}
Harris, J., Morrison, I.,
Moduli of curves. Graduate Texts in Mathematics, 187. Springer-Verlag, New York, 1998. xiv+366 pp. ISBN: 0-387-98438-0; 0-387-98429-1, MR1631825

\bibitem[K]{kass} 
J. L. Kass, 
Singular curves and their compactified Jacobians. A celebration of algebraic geometry, 391--427, Clay Math. Proc., 18, Amer. Math. Soc., Providence, RI, 2013, MR3114949

\bibitem[L]{laksov}
D. Laksov,
Weierstrass points on curves. Young tableaux and Schur functors in algebra and geometry (Toruń, 1980), pp. 221--247, Astérisque, 87–88, Soc. Math. France, Paris, 1981, MR0646822

\bibitem[LT]{laksov-thorup}
D. Laksov, A. Thorup, 
Weierstrass points and gap sequences for families of curves. Ark. Mat. 32 (1994), no. 2, 393--422, MR1318539

\bibitem[L1]{lax}
R. F. Lax,
On the distribution of Weierstrass points on singular curves. Israel J. Math. 57 (1987), no. 1, 107--115, MR0882250

\bibitem[L2]{lax2}
R. F. Lax,
Weierstrass weight and degenerations. Proc. Amer. Math. Soc. 101 (1987), no. 1, 8--10, MR0897062

\bibitem[LW]{widland-lax}
R. F. Lax, C. Widland,
Weierstrass points on Gorenstein curves. Pacific J. Math. 142 (1990), no. 1, 197--208, MR1038736

\bibitem[M1]{mumford1}
D. Mumford,
Curves and their Jacobians. The University of Michigan Press, Ann Arbor, Mich., 1975. {\rm vi}+104 pp, MR0419430

\bibitem[M2]{mumford2} 
D. Mumford,
Tata lectures on theta. II. Jacobian theta functions and differential equations. With the collaboration of C. Musili, M. Nori, E. Previato, M. Stillman and H. Umemura. Progress in Mathematics, 43. Birkhäuser Boston, Inc., Boston, MA, 1984. {\rm xiv}+272 pp. ISBN: 0-8176-3110-0, MR0742776

\bibitem[N]{neeman}
A. Neeman,
The distribution of Weierstrass points on a compact Riemann surface. Ann. of Math. (2) 120 (1984), no. 2, 317--328, MR0763909

\bibitem[Og]{og}
R. Ogawa,
On the points of Weierstrass in dimensions greater than one. Trans. Amer. Math. Soc. 184 (1973), 401--417, MR0325997 

\bibitem[Ol]{olsen}
B. A. Olsen,
On higher order Weierstrass points. Ann. of Math. (2) 95 (1972), 357--364, MR0294343

\bibitem[R]{rudnick}
Z. Rudnick,
On the asymptotic distribution of zeros of modular forms. Int. Math. Res. Not. 2005, no. 34, 2059--2074, MR2181743

\bibitem[Sh]{shivaprasad}
S. Shivaprasad,
Convergence of Bergman measures towards the Zhang measure, arXiv 2005.05753.

\bibitem[Si]{simpson}
C. T. Simpson,
Moduli of representations of the fundamental group of a smooth projective variety. I. Inst. Hautes Études Sci. Publ. Math. No. 79 (1994), 47--129, MR1307297

\bibitem[We]{wentworth}
R. Wentworth,
The asymptotics of the Arakelov-Green's function and Faltings' delta invariant. Comm. Math. Phys. 137 (1991), no. 3, 427--459, MR1105425

\bibitem[Wi]{widland}
C. Widland,
Weierstrass points on Gorenstein curves. Pacific J. Math. 142 (1990), no. 1, 197--208, MR1038736

\bibitem[Y]{yamada}
A. Yamada,
Precise variational formulas for abelian differentials. Kodai Math. J. 3 (1980), no. 1, 114--143, MR0569537

\end{thebibliography}
\end{document}